\documentclass{amsproc}%
\usepackage{amsfonts}
\usepackage{amsmath}
\usepackage{amssymb}
\usepackage{graphicx}%
\theoremstyle{plain}

\newtheorem{corollary}{Corollary}

\newtheorem{lemma}{Lemma}

\newtheorem{theorem}{Theorem}
\numberwithin{equation}{section}
\begin{document}
\title[A lower bound for the diameter of solutions to the Ricci flow]{A lower bound for the diameter of solutions to the Ricci flow with nonzero
$H^{1}(M^{n};\mathbb{R})$}
\author{Tom Ilmanen}
\address[Tom Ilmanen]{ ETH - Z\"{u}rich}
\urladdr{http://www.math.ethz.ch/\symbol{126}ilmanen}
\author{Dan Knopf}
\address[Dan Knopf]{ University of Iowa}
\email{dknopf@math.uiowa.edu}
\urladdr{http://www.math.uiowa.edu/\symbol{126}dknopf}
\thanks{Second author partially supported by NSF grant DMS - 0202796.}
\date{Version 10.02.02}
\maketitle

\section{Introduction}

Consider the Ricci flow%
\begin{equation}
\frac{\partial}{\partial t}g=-2\operatorname{Rc}(g),\qquad0\leq t<T,
\end{equation}
of a metric $g$ on $S^{1}\times S^{2}$. Our intuition suggest that no matter
how wild the metric is, the Ricci curvature in the $S^{1}$ direction should
more or less average out to zero, so that the distance around the $S^{1}$
should tend not to decrease.

In this paper, we substantiate this idea by proving the following more general
theorem. Given a Riemannian manifold $(M^{n},g)$ and a homology element
$\alpha\in H_{1}(M^{n};\mathbb{Z)}$, let $L_{\alpha}(g)$ denote the infimum of
the lengths measured with respect to $g$ of all curves representing $\alpha$.

\begin{theorem}
\label{Main}If $(M^{n},g(t):0\leq t<T)$ is a compact solution of the Ricci
flow and $\alpha\in H_{1}(M^{n};\mathbb{Z)}$ is an element of infinite order,
there exists $c=c(\alpha,g(0))\ $such that%
\[
L_{\alpha}(g(t))\geq c>0
\]
for all $t\in\lbrack0,T)$.
\end{theorem}

A particular consequence is that the diameter of $(M^{n},g(t))$ is bounded
from below independently of $t$. As a result, we can resolve a conjecture made
by Hamilton in ~\S 26 of \cite{H-95a}. Suppose that $(M^{n},g(t):0\leq t<T)$
is a solution of the Ricci flow on a maximal time interval. For $x_{j}\in
M^{n}$, $t_{j}\in\lbrack0,T)$, and $\lambda_{j}>0$, define the dilations%
\begin{equation}
g_{j}(t):=\lambda_{j}g(t_{j}+\frac{t}{\lambda_{j}}),\qquad-\lambda_{j}%
t_{j}\leq t<\lambda_{j}(T-t_{j}). \label{blowups}%
\end{equation}
If $t_{j}\nearrow T$ as $j\rightarrow\infty$ and $(M^{n},g_{j}(t),x_{j})$
converges locally smoothly to a limit $(M_{\infty}^{n},g_{\infty}%
(t),x_{\infty})$, we call the latter a \emph{final time limit flow }of
$(M^{n},g(t):0\leq t<T)$. The following result answers Hamilton's conjecture affirmatively.

\begin{corollary}
\label{Hamilton}$(S^{1}\times S^{n-1},\bar{g}(t))$ cannot arise as a final
time limit flow.
\end{corollary}

Here $\bar{g}(t)$ is the Ricci soliton%
\begin{equation}
\bar{g}(t):=ds\otimes ds+2(n-1)(\bar{T}-t)g_{\operatorname*{can}},
\label{soliton}%
\end{equation}
where $g_{\operatorname*{can}}$ is a round metric on $S^{n-1}$. Note that the
only possible final time limit flow of $(S^{1}\times S^{n-1},\bar{g}(t))$ is
$(\mathbb{R}\times S^{n-1},\bar{g}(t))$. The main content of the corollary
concerns the case that only a subsequence is known to converge.

In order to put our results into context, recall that we should not be
surprised if a solution $(M^{n},g(t))$ of the Ricci flow starting from an
arbitrary Riemannian manifold encounters a finite time singularity. Indeed,
this must be the case if the scalar curvature ever becomes everywhere
positive. (See the proof of Corollary \ref{Hamilton}, below.) To study a
finite time singularity, it is often useful to construct a sequence of
dilations (\ref{blowups}), sometimes called a blowup sequence. In certain
cases (namely if the $\lambda_{j}$ are comparable to the suprema of the
curvatures and if an injectivity radius estimate is available for the sequence
$g_{j}$) one can apply Gromov-type compactness arguments such as those in
\cite{H-95b} in order to show $C^{\infty}$ convergence to a final time limit
flow. (Such flows are also called singularity models in the literature.) Final
time limit flows have special properties which aid analysis of the original
solution $(M^{n},g(t))$. With sufficient knowledge of the limit, one can draw
useful conclusions about the analytic, geometric, and topological character of
a singular solution just prior to the formation of the singularity. The
analysis of singularities via the formation of final time limit flows is an
integral part of Hamilton's well-developed program to resolve Thurston's
Geometrization Conjecture for closed $3$-manifolds \cite{T-82} by means of the
Ricci flow. (See for example \cite{H-95a} and the survey \cite{CaoChow}.)

Theorem \ref{Main} is essentially a monotonicity result. We shall offer two
proofs which are dual to one another. The first proof (Section
\ref{Cohomology}) uses cohomology, is simpler, and is more direct. The second
proof (Section \ref{Homology}) uses homology. Its value lies in better
revealing the geometry; in particular; we hope that it is instructive in
showing how the ideas introduced here might be generalized. We briefly discuss
such potential applications in Section \ref{Conclusion}.

\section{The cohomology proof\label{Cohomology}}

Our starting point is the following observation, which was pointed out to the
first author by Sun-Chin (Michael) Chu. Let $(M^{n},g(t))$ be a solution of
the Ricci flow, and let $\phi(t)$ be a $1$-parameter family of $1$-forms
evolving by%
\[
\frac{\partial}{\partial t}\phi=\Delta_{d}\phi,
\]
where $-\Delta_{d}:=d\delta+\delta d$ is the Hodge--de Rham Laplacian.
Recalling that%
\[
\Delta_{d}\phi_{i}=\Delta\phi_{i}-R_{i}^{j}\phi_{j},
\]
where $\Delta$ is the rough Laplacian, one computes that%
\begin{align*}
\frac{\partial}{\partial t}\left\vert \phi\right\vert ^{2}  &  =\frac
{\partial}{\partial t}(g^{ij}\phi_{i}\phi_{j})\\
&  =2R^{ij}\phi_{i}\phi_{j}+2\phi^{i}\frac{\partial}{\partial t}\phi_{i}\\
&  =2\phi^{i}\Delta\phi_{i}\\
&  =\Delta\left\vert \phi\right\vert ^{2}-2\left\vert \nabla\phi\right\vert
^{2}.
\end{align*}
Applying the parabolic maximum principle, one concludes that
\begin{equation}
\left\Vert \phi(t)\right\Vert _{g(t)}\leq\left\Vert \phi(0)\right\Vert _{g(0)}
\label{key}%
\end{equation}
for as long as the solution $g(t)$ exists, where $\left\Vert \phi\right\Vert
_{g}$ denotes the supremum norm $\left\Vert \phi\right\Vert _{g}:=\sup_{x\in
M}\left\vert \phi(x)\right\vert _{g(x)}$.

We use this observation to establish a key monotonicity property. Given a
Riemannian manifold $(M^{n},g)$ and an element $\Phi$ of the first de Rham
cohomology group $H_{dR}^{1}(M^{n};\mathbb{R})$, define%
\begin{equation}
N_{g}(\Phi):=\inf_{\phi\in\Phi}\left\Vert \phi\right\Vert _{g}.
\end{equation}

\begin{lemma}
\label{Monotone1}If $(M^{n},g(t):0\leq t<T)$ is a solution of the Ricci flow,
$N_{g(t)}(\Phi)$ is a non-increasing function of time.
\end{lemma}

\begin{proof}
For any $\varepsilon>0$, there is a smooth representative $\phi_{0}\in\Phi$
such that
\[
\left\Vert \phi_{0}\right\Vert \leq N_{g(0)}(\Phi)+\varepsilon.
\]
Define $\phi(t)$ by
\begin{subequations}
\label{Hodge}%
\begin{align*}
\frac{\partial}{\partial t}\phi &  =\Delta_{d}\phi\\
\phi(0)  &  =\phi_{0},
\end{align*}
noting that a solution $\phi(t)$ exists for as long as $g(t)$ exists. Note too
that $\phi(t)\in\Phi$. Indeed, if we define a smooth function $F(t)$ by%
\end{subequations}
\begin{align*}
\frac{\partial F}{\partial t}  &  =\Delta F-\delta\phi_{0}\\
F(0)  &  =0,
\end{align*}
we have
\[
\phi=\phi_{0}+dF
\]
for all $t\in\lbrack0,T)$, because%
\[
\frac{\partial}{\partial t}(\phi_{0}+dF)=\Delta_{d}(\phi_{0}+dF).
\]
Hence by (\ref{key}), we obtain%
\[
N_{g(t)}(\Phi)\leq\left\Vert \phi\right\Vert _{g(t)}\leq\left\Vert
\phi\right\Vert _{g(0)}\leq N_{g(0)}(\Phi)+\varepsilon.
\]

\end{proof}

The following observation is of independent interest.

\begin{lemma}
\label{Norm}If $(M^{n},g)$ is a compact Riemannian manifold, then $N_{g}$ is a
norm on $H_{dR}^{1}(M^{n};\mathbb{R})$.
\end{lemma}

\begin{proof}
Homogeneity and the triangle inequality are readily verified. To show
positivity, suppose that $N_{g}(\Phi)=0$. Then there is a sequence $\left\{
\phi_{j}:j\in\mathbb{N}\right\}  \subset\Phi$ of smooth $1$-forms such that
$\left\Vert \phi_{j}\right\Vert _{g}\rightarrow0$. Fix any $\phi\in\Phi$. We
may write%
\[
\phi-\phi_{j}=dF_{j},
\]
where each $F_{j}$ is smooth. Since $M^{n}$ is compact, $\sup_{j\in\mathbb{N}%
}\left\Vert dF_{j}\right\Vert <\infty$. So after adding a locally constant
function to $F_{j}$, we may by Arzela--Ascoli select a subsequence $F_{j_{k}}$
that converges uniformly to a Lipschitz function $F$. Then%
\[
\operatorname*{ess}\sup\left\vert \phi-dF\right\vert _{g}\leq\limsup
_{k\rightarrow\infty}\left\Vert \phi-dF_{j_{k}}\right\Vert _{g}=\limsup
_{k\rightarrow\infty}\left\Vert \phi_{j_{k}}\right\Vert _{g}=0.
\]
So $dF=\phi$ almost everywhere, which implies in particular that $F$ is
smooth. Hence $\Phi=0$.
\end{proof}

We can now obtain a lower bound for the diameter of a solution of the Ricci
flow on a compact manifold $M^{n}$ with $H^{1}(M^{n};\mathbb{R})\neq\left\{
0\right\}  $.

\begin{proof}
[First proof of Theorem \ref{Main}]Let $(M^{n},g(t):0\leq t<T)$ be a solution
of the Ricci flow. Consider the natural map $\rho:H_{1}(M^{n};\mathbb{Z}%
)\rightarrow H_{1}(M^{n};\mathbb{R})$ and note that $\beta\in H_{1}%
(M^{n};\mathbb{Z})$ is a torsion element if and only if $\rho(\beta)=0$, hence
if and only if $\left\langle \Psi,\beta\right\rangle =0$ for all $\Psi\in
H^{1}(M^{n};\mathbb{R})$. So if $\alpha\in H_{1}(M^{n};\mathbb{Z})$ is an
element of infinite order, then there exists $\Phi\in H^{1}(M^{n};\mathbb{R})$
such that $\left\langle \Phi,\alpha\right\rangle >0$.

Fix any $t\in\lbrack0,T)$, and let $a$ be any curve representing $\alpha$.
Then for all $\phi\in\Phi$, we have%
\[
0<\left\langle \Phi,\alpha\right\rangle =\int_{a}\phi\leq\left\Vert
\phi\right\Vert _{g(t)}\cdot\operatorname{length}_{g(t)}(a).
\]
Taking the infimum over $\phi\in\Phi$, we get
\[
\left\langle \Phi,\alpha\right\rangle \leq N_{g(t)}(\Phi)\cdot
\operatorname{length}_{g(t)}(a)\leq N_{g(0)}(\Phi)\cdot\operatorname{length}%
_{g(t)}(a)
\]
by Lemma \ref{Monotone1}. Taking the infimum over all $a\in\alpha$, we obtain%
\begin{equation}
L_{\alpha}(g(t))\geq\frac{\left\langle \Phi,\alpha\right\rangle }%
{N_{g(0)}(\Phi)}>0. \label{lower-bound}%
\end{equation}

\end{proof}

\begin{proof}
[Proof of Corollary \ref{Hamilton}]Let $(M^{n},g(t):0\leq t<T\leq\infty)$ be a
solution of the Ricci flow on a maximal time interval, and let $g_{j}%
(t)=\lambda_{j}g(t_{j}+t/\lambda_{j})$ be a sequence of dilations such that%
\begin{equation}
(M^{n},g_{j}(t))\rightarrow(S^{1}\times S^{n-1},\bar{g}(t)),
\label{smooth-convergence}%
\end{equation}
where $\bar{g}(t)$ is defined by (\ref{soliton}). Then there exists $j_{0}$
such that $g(t_{j_{0}})$ has positive scalar curvature $R>0$. Because%
\[
\frac{\partial}{\partial t}R=\Delta R+2\left\vert \operatorname{Rc}\right\vert
^{2}\geq\Delta R+\frac{2}{n}R^{2},
\]
the maximum principle implies that the solution must fail to exist at a finite
time $T<\infty$. By Theorem 8.1 of \cite{H-95a}, a finite time singularity
implies that%
\[
\limsup_{t\nearrow T}\left(  \sup_{x\in M}\left\vert \operatorname{Rm}%
(x,t)\right\vert \right)  =\infty.
\]
Then because there is $C=C(n)$ such that
\[
\frac{\partial}{\partial t}\left\vert \operatorname{Rm}\right\vert ^{2}%
\leq\Delta\left\vert \operatorname{Rm}\right\vert ^{2}+C\left\vert
\operatorname{Rm}\right\vert ^{3},
\]
the maximum principle further implies a lower bound for the curvature blowup
rate,%
\[
\sup_{x\in M}\left\vert \operatorname{Rm}(x,t)\right\vert \geq\frac{2/C}%
{T-t}.
\]
But then smooth convergence (\ref{smooth-convergence}) is possible only if%
\[
\lim_{j\rightarrow\infty}\lambda_{j}=\infty.
\]
On the other hand, since $S^{1}\times S^{n-1}$ is compact,
(\ref{smooth-convergence}) also implies that%
\[
H_{1}(M^{n};\mathbb{Z})\cong H_{1}(S^{1}\times S^{n-1};\mathbb{Z}%
)\cong\mathbb{Z}.
\]
Let $\alpha$ generate $H_{1}(M^{n};\mathbb{Z})$. By Theorem \ref{Main}, we
have%
\[
L_{\alpha}(g(t))\geq c>0.
\]
Hence%
\[
L_{\alpha}(g_{j}(0))\geq\lambda_{j}c\rightarrow\infty
\]
as $j\rightarrow\infty$. This contradicts (\ref{smooth-convergence}) and
establishes Corollary \ref{Hamilton}.
\end{proof}

\section{The homology proof\label{Homology}}

We now seek a monotone quantity dual to the metric norms $N_{g(t)}$ defined
above on $H_{dR}^{1}$. Let $(M^{n},g)$ be a Riemannian manifold. For each free
homotopy class $\Gamma\in\operatorname*{Free}(M^{n})$, define%
\begin{align*}
\ell_{g}(\Gamma)  &  :=\inf_{\gamma\in\Gamma}\operatorname{length}_{g}%
(\gamma),\\
m_{g}(\Gamma)  &  :=\liminf_{k\rightarrow\infty}\frac{\ell_{g}(k\Gamma)}{k},
\end{align*}
where $k\Gamma$ denotes the $k$-fold cover of $\Gamma$.

We first obtain a lower bound on the decay of $\ell_{g(t)}(\Gamma)$ during the
Ricci flow.

\begin{lemma}
\label{DecayBound}Let $(M^{n},g(t):0\leq t<T)$ be a solution of the Ricci flow
and $\Gamma\in\operatorname*{Free}(M^{n})$ a free homotopy class. Then there
exists $C>0$ depending only on $n$ such that%
\[
(\ell_{g(t)}(\Gamma))^{2}\geq(\ell_{g(0)}(\Gamma))^{2}-Ct
\]
for all $t\in\lbrack0,T)$.
\end{lemma}

\begin{proof}
We may assume $\Gamma$ is nontrivial. Fix $t\in\lbrack0,T)$. There is a
nontrivial smooth closed geodesic $\gamma\in\Gamma$ such that%
\[
\operatorname{length}_{g(t)}(\gamma)=\ell_{g(t)}(\Gamma)>0.
\]
Let $V$ denote the unit tangent vector field along $\gamma$. Stability implies
that%
\begin{equation}
\int_{\gamma}(\left\vert \nabla_{V}X\right\vert ^{2}-\left\langle
R(V,X)X,V\right\rangle )\,ds\geq0 \label{Stability}%
\end{equation}
for any smooth vector field $X$ along $\gamma$. Because of holonomy, there may
not exist a parallel orthonormal frame along $\gamma$; but we can choose an
orthonormal frame $(e_{1},\ldots,e_{n})$ along $\gamma$ such that $e_{n}=V$
and
\[
\left\vert \nabla_{V}e_{i}\right\vert \leq\frac{C_{n}}{\operatorname{length}%
_{g(t)}(\gamma_{t})}=\frac{C_{n}}{\ell_{g(t)}(\Gamma)}%
\]
for $1\leq i\leq n-1$, where $C_{n}>0$ depends only on $n$. Taking $X=e_{i}$
in (\ref{Stability}) and summing over $i=1,\dots,n-1$ yields%
\[
0\leq(n-1)\left(  \frac{C_{n}}{\operatorname{length}_{g(t)}(\gamma)}\right)
^{2}\cdot\operatorname{length}_{g(t)}(\gamma)-\int_{\gamma}\operatorname{Rc}%
(V,V)\,ds.
\]
Thus%
\begin{equation}
\left.  \frac{d}{ds}(\operatorname{length}_{g(s)}(\gamma))\right\vert
_{s=t}=-\int_{\gamma}\operatorname{Rc}(V,V)\,ds\geq-\frac{(n-1)C_{n}^{2}%
}{\operatorname{length}_{g(t)}(\gamma)}. \label{bound-on-length-decrease}%
\end{equation}

Now define%
\begin{align*}
f  &  :\Gamma\times\lbrack0,T)\rightarrow\mathbb{R},\\
f(\beta,t)  &  :=\operatorname{length}_{g(t)}(\beta).
\end{align*}
Note that $f$ is continuous in $(\beta,t)$ and is $C^{1}$ in $t$ for each
fixed $\beta\in\Gamma$. Moreover, for each $u<T$, there is a compact set
$K_{u}\subseteq\Gamma$ such that%
\[
F(T):=\min_{\beta\in\Gamma}f(\beta,t)\equiv\ell_{g(t)}(\Gamma)
\]
is attained in $K_{u}$ for $0\leq t\leq u$. It follows therefore from
(\ref{bound-on-length-decrease}) that the lower derivate%
\[
\bar{D}F(t):=\liminf_{s\rightarrow t}\frac{F(s)-F(t)}{s-t}%
\]
satisfies%
\[
\bar{D}F(t)\geq-\frac{(n-1)C_{n}^{2}}{F(t)},\qquad0\leq t<T.
\]
Hence as in \S 3 of \cite{H-86}, we conclude that%
\[
(F(t))^{2}+(n-1)C_{n}^{2}t
\]
is nondecreasing, as required.
\end{proof}

The preceding lemma yields a monotonicity result dual to Lemma \ref{Monotone1}.

\begin{lemma}
\label{Monotone2}If $(M^{n},g(t):0\leq t<T)$ is a solution of the Ricci flow,
then $m_{g(t)}(\Gamma)$ is non-decreasing.
\end{lemma}

\begin{proof}
By Lemma \ref{DecayBound}, we have%
\[
(m_{g(t)}(\Gamma))^{2}=\liminf_{k\rightarrow\infty}\frac{(\ell_{g(t)}%
(k\Gamma))^{2}}{k^{2}}\geq\liminf_{k\rightarrow\infty}\frac{(\ell
_{g(s)}(k\Gamma))^{2}-Ct}{k^{2}}=(m_{g(s)}(\Gamma))^{2}%
\]
whenever $0\leq s\leq t<T$.
\end{proof}

To exploit Lemma \ref{Monotone2}, we need to know when $m_{g(0)}(\Gamma)$ is
nonzero. Let $\eta(\Gamma)\ $denote the image of $\Gamma$ in $H_{1}%
(M^{n};\mathbb{R}).$

\begin{lemma}
\label{Nontorsion}If $(M^{n},g)$ is a Riemannian manifold and $\Gamma
\in\operatorname*{Free}(M^{n})$ is a free homotopy class such that
$\eta(\Gamma)$ is nonzero, then $m_{g}(\Gamma)>0$.
\end{lemma}

\begin{proof}
Since $\eta(\Gamma)\neq0$, there exists $\Phi\in H^{1}(M^{n};\mathbb{R})$ such
that $\left\langle \Phi,\eta(\Gamma)\right\rangle >0$. For any $\phi\in\Phi$
and any curve $\gamma\in k\Gamma$, we have%
\[
\left\langle \Phi,\eta(k\Gamma)\right\rangle =\int_{\gamma}\phi\leq\left\Vert
\phi\right\Vert \cdot\operatorname{length}_{g}(\gamma).
\]
Taking the infimum over $\phi$ and $\gamma$ yields%
\[
\left\langle \Phi,\eta(k\Gamma)\right\rangle \leq N_{g}(\Phi)\cdot\ell
_{g}(k\Gamma).
\]
Hence%
\[
m_{g}(\Gamma)=\liminf_{k\rightarrow\infty}\frac{\ell_{g}(k\Gamma)}{k}%
\geq\liminf_{k\rightarrow\infty}\frac{\left\langle \Phi,\eta(k\Gamma
)\right\rangle }{kN_{g}(\Phi)}=\frac{\left\langle \Phi,\eta(\Gamma
)\right\rangle }{N_{g}(\Phi)}>0.
\]

\end{proof}

These observations lead to another proof of the main result of this paper.

\begin{proof}
[Second proof of Theorem \ref{Main}]Let $(M^{n},g(t):0\leq t<T)$ be a solution
of the Ricci flow, and let $\alpha\in H_{1}(M^{n};\mathbb{Z})$ be an element
of infinite order. Then there exists a free homotopy class $\Gamma
\in\operatorname*{Free}(M^{n})$ whose image in $H_{1}(M^{n};\mathbb{Z})$ is
$\alpha$. Clearly, $L_{\alpha}(g(t))=\ell_{g(t)}(\Gamma)\geq m_{g(t)}(\Gamma
)$. Since $\alpha$ is of infinite order, $\eta(\Gamma)\in H_{1}(M^{n}%
;\mathbb{R})$ is nonzero. So we can apply Lemmas \ref{Monotone2} and
\ref{Nontorsion} to conclude that%
\[
L_{\alpha}(g(t))\geq m_{g(t)}(\Gamma)\geq m_{g(0)}(\Gamma)>0.
\]

\end{proof}

\section{Concluding remarks\label{Conclusion}}

Although one expects Ricci flow evolutions to encounter finite-time
singularities for a large class of initial Riemannian manifolds, the main
result of this paper shows that there are topological restrictions on the
geometry of such singularities. Motivated by this observation, we pose the
following problems.

\medskip

\noindent\textbf{Problem 1. }Suppose that $(M^{n},g_{j}(t))$ is a blowup
sequence converging smoothly (in the pointed category) to a solution
$(M_{\infty}^{n},g_{\infty}(t))$ of the Ricci flow. Show that the image of
$H^{1}(M_{\infty}^{n};\mathbb{Z})$ in $H^{1}(M^{n};\mathbb{Z})$ under the
natural map is finite.

\medskip

\noindent\textbf{Problem 2. }The lens spaces $L(p,q)$ demonstrate that there
can be no lower bound for the length of a torsion element, hence no torsion
analogue of Theorem \ref{Main}. If there is torsion in $H_{1}(M^{n}%
;\mathbb{Z})$, is it true that any solution $(M^{n},g(t))$ of the Ricci flow
must become singular in finite time?

\end{document}